\documentclass[12pt]{amsart}

\newcommand{\Q}{\mathbb Q}
\newcommand{\R}{\mathbb R}
\newcommand{\Z}{\mathbb Z}

\renewcommand{\epsilon}{\varepsilon}
\renewcommand{\phi}{\varphi}

\newtheorem{Theorem}{Theorem}[section]

\newtheorem{Definition}{Definition}[section]

\newtheorem{Remark}{Remark}[section]

\begin{document}
\title[A.Borisov, Toric Singularities]{ On Classification of Toric Singularities  }
\author[A.Borisov, Toric Singularities]{  Alexandr Borisov}
\address{Department of Mathematics, Penn State University, University Park, PA $16802$, USA
}
\email{borisov@math.psu.edu}
\date{Apr 12, 1998}
\maketitle
\thispagestyle{empty}












\section{Introduction}

In 1988 S. Mori, D. Morrison, and I. Morrison published a paper  \cite{MMM} in which they gave a computer-based  conjectural classification of four-dimensional terminal cyclic quotient singularities of prime index. This classification was partially proven in 1990 by G. Sankaran (cf. \cite{Sankaran}). Basically, the conjecture of \cite{MMM}  says that these singularities form a finite number of series (cf. definition in section 2 below) with the precise list of series given.

While a precise list of series in higher dimensions would be way too long, the qualitative version of this conjecture makes sense for toric singularities of any dimension. It turned out, to my great surprise, that this and much more was basically proven in 1991 in a beautiful paper by Jim Lawrence (cf. \cite{Lawrence}). His proof and motivation came from the geometry of numbers and he was obviously unaware of the papers of Mori-Morrison-Morrison and Sankaran.

So in this short note I just bring this all together. In section 2 the necessary definitions are given and the main theorem is proven. In section 3 some interesting related open problems are discussed.

{\bf Acknowledgments.} I can claim only little credit for the proof of the main theorem, as this is just an algebro-geometric corollary of the result of Jim Lawrence. I discovered the paper of Lawrence using MathSciNet.

\section{Series of toric singularities}

First of all let me recall (cf. \cite{D}) that an $n-$dimensional $\Q-$factorial toric singularity is given by a rational simplicial cone in $\R^n$. If we identify by a linear transformation the closest integral points of the one-dimensional faces of this cone with the standard basis $(e_i)$ of $\R^n,$ the original lattice of integral points $L$ will contain $\Z^n$ as a sublattice of finite index. So $L/ \Z^n $ is a finite subgroup of $T^n=\R^n/ \Z^n .$ In the opposite direction, any finite subgroup $\bar{L} \subset T^n$ gives rise to some lattice $L.$ And if additionally $\R \cdot e_i \cap L =\R \cdot e_i \cap \Z^n$  for every $i$ we get a toric singularity using the standard cone in $\R^n$ and the lattice $L.$


Let me introduce the following definition.

\begin{Definition}
A series of $\Q-$factorial toric singularities is a closed subgroup $V \subset T^n$ together with finitely many proper closed subgroups $V_i\subset V$ and $V_{ij}\subset V_i.$ A toric singularity belongs to this series if and only if it is defined by a finite subgroup $G\subset V$ such that
$G\cap V_i \subset \bigcup \limits_j V_{ij}$
for all $i$.
\end{Definition}

\begin{Remark} I use  the word ``series" instead of the word ``family" to avoid confusion with the standard use of the word  ``family" in algebraic geometry.
\end{Remark}

\begin{Remark}
If we restrict our attention to the cyclic quotient singularities of prime index we can simplify the above condition considerably. Namely, if $P\in T^n$ is a generator of the corresponding finite subgroup, the above condition could be replaced by 
$$P\in V \setminus  \bigcup \limits_i \Big( V_i \setminus \bigcup \limits_j V_{ij}\Big).$$
\end{Remark}
This agrees nicely with the conjecture of Mori-Morrison-Morrison. However in general we do need a more complicated definition as above.

\begin{Definition}
The dimension of the series as above is $\dim V.$
\end{Definition}

\begin{Remark}
In the terminology of \cite{MMM} the terminal cyclic quotient singularity of prime index is stable if it belongs to a family of terminal toric singularities of dimension at least one.
\end{Remark}

It was conjectured in \cite{MMM} that in dimension four all but finitely many cyclic quotient singularities of prime index are stable with a precise finite list of stable series and exceptions. The list of stable series, i.e. the series of dimension at least one, was proven to be complete by G. Sankaran (cf. \cite{Sankaran}). The list of exceptions, i.e. zero-dimensional series is too long to be proven complete without an extensive computer calculation. Even with a computer it is not exactly clear how to do it. As far as I know it is not yet established.

The main result of this paper is the following.
\begin{Theorem}
For every $\epsilon > 0$ consider all $\Q-$factorial toric singularities with minimal log-discrepancy  greater than (greater than or equal to) $\epsilon .$ Then they form a finite number of series.
\end{Theorem}

{\bf Proof.} First of all, let's recall that the minimal log-discrepancy of a toric singularity is equal to the minimal sum of coordinates of the non-zero elements in the corresponding finite subgroup of $T^n,$ where $T^n$ is identified with the set of points $(x_i)$ such that $0\le x_i <1$ (cf. \cite{Reid}, \cite{Bor?}).

For every $\epsilon >0$ let's define two subsets of $T^n$ as follows. 
 $$S_{\epsilon }= \{(x_i)\in T^n: x_i>0, \sum x_i <\epsilon \}$$
 $$S^{'}_{\epsilon }= \{(x_i)\in T^n: x_i>0, \sum x_i \le \epsilon \}$$ 
In the terminology of Lawrence (cf. \cite{Lawrence}, section 3) these subsets are full. By the result of Lawrence (cf. \cite{Lawrence}, Thm 1) there are just finitely many closed subgroups $V^{(\nu)}$ of $T^n$ which are maximal among those that avoid $S_{\epsilon }$ ($S^{'}_{\epsilon }$). So any finite subgroup that defines a toric singularity with minimal log-discrepancy at least $\epsilon $ (greater than $\epsilon $) must be contained in some $V^{(\nu)}=V.$

Unfortunately, not every finite subgroup $F$ of $V$ defines a toric singularity with the desired property. First of all we have the condition that for all $i=1,2,...,n$
$$F\cap (\R/\Z ) \cdot e_i =\{0\}$$
which is equivalent to saying that $F$ defines a toric singularity. Also $F$ could contain a point on the boundary of $T^n$ with the sum of coordinates less than (or equal to) $\epsilon .$ To deal with this, let's consider separately all lower-dimensional faces of $T^n$, containing $0$, of dimension $2,3,..., n-1.$ Let $L$ be one such face. Let's write down what it means that $F$ has a point {\it strictly inside } L with the sum of coordinates less than (or equal to) $\epsilon .$

Let's denote by $I_L \subset \{1,2,...,n\}$ the set of coordinates that are identically zero on $L.$ Then the above condition means that
$F\cap S_L =\emptyset$ ($F\cap S^{'}_L =\emptyset$) where
$$S_L= \{(x_i)\in L: x_i>0 \,\,\, \forall i\notin I_L, \sum x_i <\epsilon \}$$
$$S^{'}_L= \{(x_i)\in L: x_i>0 \,\,\, \forall i\notin I_L, \sum x_i \le \epsilon \}$$
Again by the result of Lawrence there are just finitely many closed subgroups $V_{L,\tau}$ of $V\cap L$ which are maximal among those that avoid $S_L$ ($S^{'}_L$). So the condition in question can be written as 
$$F\cap (V\cap L)\subset \bigcap \limits_{\tau} V_{L,\tau}.$$
I should note that $V$ could be itself a face of $T^n$. This happens, for example for all $\epsilon >\frac12 .$ Then formally we should be careful because $V\cap L$ is not a proper subgroup of $V$ so the above condition is a bit different from the Definition 2.1. In fact in this case we get not one family generally defined by $V$ but several smaller families. This can be easily handled by the induction on $n.$ 

With the above remark, the proof of the theorem is now complete.

\begin{Remark}One immediate corollary of this theorem, and the result of G. Sankaran (\cite{Sankaran}) is that the list of Mori-Morrison-Morrison (\cite{MMM}) is complete up to possibly a finite number of overlooked exceptions. It also gives a similar but qualitative result for the terminal cyclic quotient singularities of prime index in any dimension, which was also conjectured in \cite{MMM}.
\end{Remark}

\begin{Remark} One can also derive from the Theorem 2.1 Shokurov's conjecture for toric singularities, that their minimal log-discrepancies could only accumulate from above (cf. \cite{Lawrence}, proof of Lemma 5). This is one of the results of my paper \cite{Bor?} which I published before discovering the Lawrence's paper. However the main result of \cite{Bor?}, i.e. that the limit points are essentially the minimal log-discrepancies of cyclic quotient singularities of lower dimension, doesn't seem to follow  right away from the Theorem 2.1.
\end{Remark}

\section{Open questions}
There are many related open questions that arise. Let me list here those that I think are the most interesting.

{\bf Question 1.} When $n$ and/or $\epsilon$ varies how does the finite set of subgroups from the Theorem 2.1 vary? What is its combinatorial and geometric meaning?
Even when $n=2$ the question is interesting and in general it is wide open. One can define the Hilbert polynomial $H$ of this constructive set by using the polynomial measure which is defined by the property $H(V,x)=r\cdot x^{\dim V}$ if $V$ is a subgroup of $T^n$ of dimension $\dim V$ with $r$ connected components. It could be interesting to calculate this polynomial in certain cases.

{\bf Question 2.} It would be of some interest to write computer codes that would effectively find this finite set of series for different $n$ and $\epsilon.$ This question was already asked by Lawrence in \cite{Lawrence}. Some computer-free results in this direction are contained in \cite{MMM}, \cite{Sankaran}, \cite{MorStev}.

{\bf Question 3.} One of the results of Lawrence is the following.

{\bf Theorem} (\cite{Lawrence}, Corollary 1.B) {\it Let $S$ be a closed subset of $T^n$ such that if $x\in S$  and $m$ is a positive integer then $mx\in S.$ Then $S$ is a finite union of closed subgroups of $T^n.$ }

If one tries to generalize this to noncommutative compact Lie groups instead of $T^n$ it fails in general. For example, we can take $S$ to be a (closed subset of) the set of all operators $A$ in $SO_n$ with the property $A^2=Id.$ It would be interesting to determine exactly the structure of such sets for all compact Lie groups.

{\bf Question 4.} There is a well-known classification of three-dimensional terminal singularities, due to S. Mori and  Miles Reid (cf. \cite{Reid}). It is complete modulo some unanswered questions about terminality of specializations of some exceptional families. Basically the answer is that these singularities form a finite number of ``series". I wonder if the something like this is true in higher dimensions, and if one replaces terminal by $\epsilon-$logterminal singularities. It is of course not clear at all what to mean by a ``series" in this general setting.

\end{document}